\documentclass[11pt]{article}

 \usepackage{setspace}
\usepackage{natbib}

\usepackage{latexsym}
\usepackage{calc}

\usepackage{amssymb}
\usepackage{amsmath}
\usepackage{graphicx}

\usepackage{psfrag}

\usepackage{float}

\usepackage{color}

\usepackage{nicefrac}

\usepackage{ntheorem}

\usepackage{hyperref}

\newcounter{hours}\newcounter{minutes}

\def\nr{\par \noindent}

\def\beq{\begin{equation}}
\def\eeq{\end{equation}}

\newtheorem{theorem}{Theorem}
\newtheorem*{btheorem}{Blackwell's Theorem:}
\newtheorem{lemma}{Lemma}
\newtheorem{corollary}{Corollary}

\newtheorem{proposition}{Proposition}
\newtheorem{assumption}{Assumption}
\newtheorem{definition}{Definition}

\newtheorem{example}{Example}
\newtheorem{remark}{Remark}
\newcommand{\proof}{\bf Proof: \rm \nr}

\newcommand{\qed}{\hfill $\Box$ \nr \medskip}

\def\ba{\begin{array}}
\def\ea{\end{array}}
\def\beann{\begin{eqnarray*}}
\def\eeann{\end{eqnarray*}}
\def\bea{\begin{eqnarray}}
\def\eea{\end{eqnarray}}

\def\BT{\begin{theorem}}
\def\ET{\end{theorem}}
\def\BL{\begin{lemma}}
\def\EL{\end{lemma}}
\def\BC{\begin{corollary}}
\def\EC{\end{corollary}}
\def\BE{\begin{example}}
\def\EE{\end{example}}
\def\BD{\begin{definition}}
\def\ED{\end{definition}}
\def\BR{\begin{remark}}
\def\ER{\end{remark}}
\def\BAS{\begin{assumption}}
\def\EAS{\end{assumption}}
\def\BI{\begin{itemize}}
\def\EI{\end{itemize}}

\def\BMP{\begin{minipage}{9.5cm}}
\def\EMP{\end{minipage}}
\def\MPT{\begin{minipage}{11.5cm}}
\def\EPT{\end{minipage}}

\def\R{\mathbb{R}}

\newcommand*\tr[1]{\mathop{}\!\mathrm{tr}\left(#1\right)}

\newcommand*\samethanks[1][\value{footnote}]{\footnotemark[#1]}

\title{Shannon's comparison of channels \\ characterized by optimal decision making}
\author{S. L\"ammel
\thanks{
Department of Mathematics, Chemnitz University of Technology,
Reichenhainer Str. 41, 09126
Chemnitz, Germany; e-mail: sebastian.laemmel@mathematik.tu-chemnitz.de (corresponding author), vladimir.shikhman@mathematik.tu-chemnitz.de.
 } \and V. Shikhman\samethanks[1]
}

\begin{document}
\maketitle
\vspace{-5ex}
\abstract{ 
According to Blackwell's Theorem it is equivalent to compare channels by either a garbling order or optimal decision making.
This equivalence does not hold anymore if also allowing pre-garbling, i.\,e. for the so-called Shannon-order \citep[see][]{rauh:2017}. We show that the equivalence fails in general even if the set of decision makers is reduced.
This is overcome by the introduction of convexified Shannon-usefulness as a preference relation of decision makers over channels. We prove that convexified Shannon-order and convexified Shannon-usefulness are equivalent.
}

\vspace{2ex}

\section{Introduction}
\label{sec:intro}

In order to compare channels, one has to rely on e.\,g. the channel capacity. However, the comparison of channel capacities does not fit into a decision framework. 
The latter is based on the comparison of expected utilities attached to channels by decision makers.
This approach has been introduced by Bohnenblust, Shapley and Sherman in the context of experiments \citep[see][]{lecam:1996}. 
Let a channel be represented by a column stochastic matrix $C \in \R^{m \times n}$, where
$m$ and $n$ are the lengths of the output and input alphabet, respectively.
The input distribution is stored in the diagonal matrix $\mathit{\Pi} \in \R^{n \times n}$. 
Suppose a decision maker chooses an action based on the output of the channel in order to maximize utility. For that, we call a function, which associates the channel's output with the set of feasible actions, a strategy. 
Any strategy can be represented without loss of generality by a stochastic matrix $A \in \R^{m \times m}$. Indeed, the number of actions can be enlarged to the length of the output alphabet $m$. The set of quadratic stochastic matrices will be denoted by $\mathcal A$. 
The set of all joint distributions of actions and inputs when endowed with $C$ is therefore represented by the so-called policy space \citep{perez-richet:2017}:
\[
\mathit{\Phi}(C) = \left\{\left.D \in \R^{m \times n}\right\vert \mbox{There exists } A \in {\mathcal A} \mbox{ with } D=A \cdot C \right\}.
\]
Furthermore, let $U \in \R^{n \times m}$ be the utility matrix of a particular decision maker. %The set of utility matrices will be denoted by $\mathcal U$. 
A decision maker maximizes utility under all possible joint distributions $D$. Therefore, the following optimization problem is to be solved:
\[
   \max_{D \in \mathit{\Phi}(C)}\tr{U\cdot D \cdot \mathit{\Pi}},
\]
where $\tr{\cdot}$ denotes the trace of a matrix.
This optimization problem provides a comparison of channels by usefulness.

\begin{definition}[Blackwell-usefulness, \citet{blackwell:1953}]
Let  $C, \bar C \in \R^{m \times n}$ be two channels with the same input alphabet. We say that $C$ is more Blackwell-useful than $\bar C$ (denoted by $C \succcurlyeq_B \bar C$) if 
for all utility matrices $U \in \R^{n \times m}$ it holds:
\[
\max_{D \in \mathit{\Phi}(C)}\tr {U\cdot D \cdot \mathit{\Pi}}\ge \max_{D \in \mathit{\Phi}\left(\bar C\right)}\tr {U\cdot D \cdot \mathit{\Pi}}.
\]
This means that every decision maker gains by using $C$ at least the utility he or she would gain by $\bar C$.
\end{definition}

Another way to compare channels focuses on the possibility of reproducing one channel by another.
\begin{definition}[Blackwell-order, \citet{blackwell:1953}]
Let $C, \bar C \in \R^{m \times n}$ be two channels with the same input alphabet. We say that $\bar C$ is a garbling of $C$ (denoted by $C \trianglerighteq_B \bar C$) if 
there exists a stochastic matrix $M \in \R^{ m \times m} $ with
\[
\bar C=M \cdot C.
\]
We call $\trianglerighteq_B$ the partial Blackwell-order of channels.
\end{definition}

 \cite{blackwell:1953} showed the equivalence of Blackwell-usefulness and Blackwell-order.
%Those apply for channels with the same input alphabet as well.\par
%As before mentioned comparison by Blackwell-usefulness and Blackwell-order are equivalent.

\begin{btheorem}
It  holds for channels $C, \bar C \in \R^{m \times n}$ with the same input alphabet:
\[
C \trianglerighteq_B \bar C \quad 
\Leftrightarrow \quad
C \succcurlyeq_B \bar C.
\]
\end{btheorem}
Apparently, not all channels are comparable. It appears reasonable to enlarge the definition of garbling by allowing pre-garbling additionally. 
Garbling corresponds then to decoding, and pre-garbling to coding.
%The decision maker can effect both, input and output of a channel.
\begin{definition}[Shannon-order,  \citet{shannon:1958}]
Let $C, \bar C \in \R^{m \times n}$ be two channels with the same input alphabet. We say $\bar C$ is a Shannon-garbling of $C$ (denoted by $C \trianglerighteq_S \bar C$) if there exist stochastic matrices $M \in \R^{m \times m}$ and $N \in \R^{n \times n}$ with
\[
\bar C= M \cdot C \cdot N.
\]
We call $\trianglerighteq_S$ the partial Shannon-order of channels.
\end{definition}
This new partial ordering is not only finer than the Blackwell-order, but it also appears more suitable for channels. \cite{shannon:1958} introduced it and interpreted $C \trianglerighteq_S \bar C$ as $\bar C$ is included in $C$.

This paper studies if there is an appropriate definition of Shannon-usefulness that provides the equivalence to the Shannon-order in a similar way as in Blackwell's Theorem.
For this purpose, we introduce in Section 2 the reduced Blackwell-usefulness which generalizes the Blackwell-usefulness and is defined with respect to subsets of utility matrices. Furthermore, we consider some particular subsets of utility matrices and examine the relation between the corresponding reduced Blackwell-usefulness and the Shannon-order. The main result of Section 2 will be Theorem \ref{th:buso}, which states that there is no equivalence between the latter for any subset of utility matrices and channels of size $\R^{m \times (2^{m-2}+1)}$ with $m\ge 2$. Section 3 introduces the notion of convexified Shannon-usefulness. Theorem \ref{th:shannonstheorem} states that convexified Shannon-usefulness and convexified Shannon-order (as introduced by \cite{shannon:1958}) are equivalent.
\section{Blackwell-usefulness and Shannon-order}
\label{sec:neg}

First, we recall that Blackwell-usefulness is not preserved by the Shannon-order as the following example by \cite{rauh:2017} shows.
\begin{example}[Failure of Blackwell-usefulness, \citet{rauh:2017}]
\label{ex:rauh}
We consider two channels
\[
C = \begin{pmatrix}\nicefrac{9}{10} & 0 \\ \nicefrac{1}{10} & 1\end{pmatrix}, 
\quad 
\bar C = \begin{pmatrix}0 & \nicefrac{9}{10} \\ 1 &  \nicefrac{1}{10}\end{pmatrix}
\]
with the same uniformly distributed input alphabet, i.\,e.
\[
\mathit{\Pi}=\begin{pmatrix} \nicefrac{1}{2} &0 \\ 0 & \nicefrac{1}{2}\end{pmatrix}.
\]
It is straightforward to see that $\bar C$ is a Shannon-garbling of $C$ with 
\[
   M=\begin{pmatrix} 1 &0 \\ 0 & 1\end{pmatrix}, \quad 
   N= \begin{pmatrix}0 & 1 \\ 1 & 0 \end{pmatrix}.
\]
However, the maximal expected utility of $\bar C$ is greater than that of $C$, at least for the following utility matrix
\[
U = \begin{pmatrix} 2 & 0 \\ 0 & 1\end{pmatrix}.
\]
It holds namely:
\[
\max_{D \in \mathit{\Phi}(C)}\tr {U\cdot D \cdot \mathit{\Pi}}= \nicefrac{28}{20}, \quad 
\max_{D \in \mathit{\Phi}\left(\bar C\right)}\tr {U\cdot D \cdot \mathit{\Pi}}=\nicefrac{29}{20}.
\]
Thus, $C$ is not more Blackwell-useful than $\bar C$.
\qed
\end{example}
Therefore, Shannon-order and Blackwell-usefulness cannot be equivalent for channels with the same input alphabet. Let us modify the definition of Blackwell-usefulness instead. As a starting point we reduce the set of utility matrices and, thus, decision makers, for which the maximal expected utility is compared. 

\begin{definition}[Reduced Blackwell-usefulness]
Let $C, \bar C \in \R^{m \times n}$ be two channels with the same input alphabet and ${\mathcal U} \subseteq \R^{n \times m}$ a subset of utility matrices.
We say that $C$ is more Blackwell-useful than $\bar C$ with respect to $ {\mathcal U}$ (denoted by $C \succcurlyeq_{B}^{\mathcal U} \bar C$) if for all utility matrices $U \in {\mathcal U}$ it holds:
\[
\max_{D \in \mathit{\Phi}(C)}\tr {U\cdot D \cdot \mathit{\Pi}}\ge \max_{D \in \mathit{\Phi}\left(\bar C\right)}\tr {U\cdot D \cdot \mathit{\Pi}}.
\]
This means that every decision maker endowed with a utility matrix $U \in {\mathcal U}$ gains by using $C$ at least the utility he or she would gain by using $\bar C$.
\end{definition}
We are aiming to identify a suitable subset ${\mathcal U}$ of utility matrices for which the reduced Blackwell-usefulness characterizes the Shannon-order, i.\,e.
\begin{equation}
\label{eq:buso}
C \trianglerighteq_S \bar C \quad 
\Leftrightarrow \quad
C \succcurlyeq_B^{{\mathcal U}} \bar C.
\end{equation}
For this purpose we assume throughout this section that the input alphabet is uniformly distributed, i.\,e.
\[
    \mathit{\Pi} = \mbox{diag}\left(\frac{1}{n}, \ldots, \frac{1}{n}\right).
\]
The case of a general $\mathit{\Pi}$ is covered in Remark \ref{re:allpi} below.

Let us define some subsets of utility matrices used in what follows.
\begin{definition}[Subsets of utility matrices]
\label{def:utility}
We call a utility matrix $U$
\begin{itemize}
\item[(1)] indifferent, if all its columns are identical;
\item[(2)] exact, if it is a positive multiple of a permutation matrix;
%i.e. there exists $a>0$ and a permutation matrix $P$, such that $U= a \cdot P$;
\item[(3)] oblivious, if it is a positive multiple of a matrix whose columns are coordinate vectors.
\end{itemize}
The sets of indifferent, exact, and oblivious utility matrices will be denoted by $\mathcal I$, $\mathcal E$, and $\mathcal O$, respectively.
Furthermore, we denote by $\mathcal D$ the set of positive multiples of doubly-stochastic matrices.
\end{definition}

Next, we justify the utility notions from Definition \ref{def:utility}.

\begin{remark}[Indifferent utility]
\label{rem:indif}
A decision maker endowed with an indifferent utility matrix $U$ will achieve the maximal expected utility independently of the channel and the chosen action, since
\[
  \tr{U \cdot D \cdot \mathit{\Pi}} = \frac{1}{n}\tr{U}
\]
is  constant for any $C$ and $D \in \Phi(C)$. This means that the change of either the channel or the action is redundant. \qed
\end{remark}

\begin{remark}[Exact utility]
In contrast, let a decision maker be endowed with an exact utility matrix $\tilde U=\alpha \cdot \tilde P$, where $\tilde P \in \mathbb{R}^{n \times n}$ is a permutation matrix. Then, it holds for $D = A \cdot C$ with $A \in \mathcal{A}$:
\[
  \tr {\tilde U\cdot D \cdot \mathit{\Pi}} = 
  \frac{1}{n} \tr {\tilde U\cdot A \cdot C}
  = \frac{1}{n} \tr {C \cdot \tilde U \cdot A}
  = \frac{\alpha}{n} \tr {\tilde C \cdot A},
\]
where the channel $\tilde C = C \cdot \tilde P$ emerges from $C$ by $\tilde P$-permutation of the input alphabet.  
Hence, we have with $C= \left(c_{ij}\right)$ and $\tilde C= \left(\tilde c_{ij}\right)$:
\[
  \max_{D \in \mathit{\Phi}(C)}\tr {\tilde U\cdot D \cdot \mathit{\Pi}} = 
  \frac{\alpha}{n}\max_{A \in \mathcal{A}}  \tr {\tilde C \cdot A} = 
  \frac{\alpha}{n} \sum_{i=1}^{n} \max_{1 \leq j \leq n} \tilde c_{ij} =
  \frac{\alpha}{n} \sum_{i=1}^{n} \max_{1 \leq j \leq n} c_{ij}. 
\]
Here, the maximum expected utility is achieved by taking e.\,g. $A= \left(a_{ij}\right)$ with $a_{ji} = 1$ for exactly one index $j$ with $\displaystyle \max_{1 \leq k \leq n } c_{ik} = c_{ij}$, and $a_{ji}=0$ otherwise.
The latter means that the decision maker chooses the action which corresponds to the 
largest transmission probability of the channel's output.\qed
\end{remark}

\begin{remark}[Oblivious utility]
Analogously, let a decision maker be endowed with an oblivious utility matrix $U' = \alpha \cdot S'$, where the columns of $S' \in \mathbb{R}^{n \times m}$ are coordinate vectors. Then, it holds for $D = A \cdot C$ with $A \in \mathcal{A}$:
\[
  \tr {U'\cdot D \cdot \mathit{\Pi}} = 
  \frac{1}{n} \tr {U'\cdot A \cdot C}
  = \frac{1}{n} \tr {C \cdot U' \cdot A}
  = \frac{\alpha}{n} \tr {C' \cdot A},
\]
where the channel $C' = C \cdot S'$ emerges from $C$ by replacing 
some letters of the input alphabet by the others, or by deleting them.  
Hence, we have with $C= \left(c_{ij}\right)$ and $ C'= \left( c'_{ij}\right)$:
\[
  \max_{D \in \mathit{\Phi}(C)}\tr {U'\cdot D \cdot \mathit{\Pi}} = 
  \frac{\alpha}{n}\max_{A \in \mathcal{A}}  \tr {C' \cdot A} = 
  \frac{\alpha}{n} \sum_{i=1}^{m} \max_{1 \leq j \leq m} c'_{ij} \leq
  \frac{\alpha}{n} \sum_{i=1}^{m} \max_{1 \leq j \leq n} c_{ij}. 
\]
Here, the maximum expected utility is achieved by taking e.\,g. $A= \left(a_{ij}\right)$ with $a_{ji} = 1$ for exactly one index $j$ with $\displaystyle \max_{1 \leq k \leq m } c'_{ik} = c'_{ij}$, and $a_{ji}=0$ otherwise.
The latter means that the decision maker chooses the action which corresponds to the 
largest remaining transmission probability of the channel's output. \qed
\end{remark}

We now focus on the sufficiency of the Shannon-order for the reduced Blackwell-usefulness.
For the set $\mathcal{I}$ of indifferent utility matrices this is trivially true due to Remark \ref{rem:indif}. It turns out that this is also true for the set $\mathcal{E}$ of exact utility matrices.

\begin{proposition}[Sufficiency for $\mathcal{E}$]
\label{prop:exact}
Let ${\mathcal U}$ be a subset of $\mathcal E$. Then, for every two channels $C, \bar C \in \R^{n\times n}$ with the same input alphabet it holds:
\[
   C \trianglerighteq_S \bar C \quad \Rightarrow \quad C \succcurlyeq_B^{{\mathcal U}} \bar C.
\]

%\[
%C \trianglerighteq_S \bar C \Rightarrow C \succcurlyeq_B^{{\mathcal P}} \bar C.
%\]
\end{proposition}
\proof
Let $U \in {\mathcal U}$. There exists $\alpha>0$ and a permutation matrix $P$ such that $U= \alpha \cdot P$. Let $\bar C = M \cdot C \cdot N$.
We define a subset of stochastic matrices 
\[
\mbox{$\mathcal B$} = \left\{ B \, \left\vert \, \mbox{There exists } A \in \mbox{$\mathcal A$ with } B=P^T \cdot N \cdot P \cdot A \cdot M \right. \right\}.
\]
It follows:
\begin{align}
\max_{D \in \mathit{\Phi}(C)}\tr {U\cdot D \cdot \mathit{\Pi}} &= \frac{\alpha}{n}\max_{A \in \mbox{$\mathcal A$}} \tr {P\cdot A \cdot C} \notag \\
&\ge  \frac{\alpha}{n} \max_{B \in \mbox{$\mathcal B$}} \tr {P\cdot B \cdot C} \notag \\
&= \frac{\alpha}{n} \max_{A \in \mbox{$\mathcal A$}} \tr {P\cdot P^T \cdot N \cdot P \cdot A \cdot M \cdot C} \notag \\
&= \frac{\alpha}{n} \max_{A \in \mbox{$\mathcal A$}} \tr {P\cdot A \cdot \bar C} \notag \\
&= \max_{D \in \mathit{\Phi}\left( \bar C \right)}\tr {U\cdot D \cdot \mathit{\Pi}}\notag. 
\end{align}\qed

We now try to extend the set $\mathcal{E}$, such that the sufficiency part remains true. One possibility is to analyze the larger set of oblivious utility matrices $\mathcal{O} \supset \mathcal{E}$. Alternatively, we examine the set of positive multiples of doubly-stochastic utility matrices $\mathcal{D} \supset \mathcal{E}$.

\begin{example}[Failure of sufficiency for $\mathcal O$ and $\mathcal D$]
The implication in (\ref{eq:buso}) does not  hold in general either for ${\mathcal U}=\mathcal O$ or for ${\mathcal U}=\mathcal D$.
\label{ex:randd}
\begin{itemize}
\item[(1)] Suppose that we have the oblivious utility matrix
\[
U = \begin{pmatrix}  1 & 1 & 0 \\  0& 0 & 1\\ 0& 0 & 0\\ \end{pmatrix} \in \mathcal O,
\]
and two channels
\[
C=\begin{pmatrix}  0 & 0 & 1 \\ \nicefrac{1}{2} & \nicefrac{1}{2} & 0\\ \nicefrac{1}{2} & \nicefrac{1}{2} & 0\\ \end{pmatrix}, \quad
\bar C = \begin{pmatrix}  0 & 0 & 1 \\  \nicefrac{1}{2} & \nicefrac{1}{2} & 0\\ \nicefrac{1}{2} & \nicefrac{1}{2} & 0\\ \end{pmatrix} 
\cdot 
\begin{pmatrix}  0 & 0 & 1 \\ 0 & 1 & 0\\ 1 & 0 & 0\\ \end{pmatrix} 
= 
\begin{pmatrix}  1 & 0 & 0 \\ 0& \nicefrac{1}{2} & \nicefrac{1}{2}\\ 0& \nicefrac{1}{2} & \nicefrac{1}{2}\\ \end{pmatrix}.
\]
Hence, $\bar C$ is a Shannon-garbling of $C$. Assuming the uniform distribution of the input alphabet, the maximal expected utilities of $C$ and $\bar C$ are
\[
\max_{D \in \mathit{\Phi}(C)}\tr {U\cdot D \cdot \mathit{\Pi}}= \nicefrac{1}{3}, \quad 
\max_{D \in \mathit{\Phi}\left(\bar C\right)}\tr {U\cdot D \cdot \mathit{\Pi}}=\nicefrac{2}{3},
\]
respectively. Thus, $C$ is not more Blackwell-useful with respect to $\mathcal O$ than $\bar C$.

\item[(2)] Suppose that we have the utility matrix
\[
U = \begin{pmatrix}  1 & 0 & 0 \\ 0& \nicefrac{1}{2} & \nicefrac{1}{2}\\ 0& \nicefrac{1}{2} & \nicefrac{1}{2}\\ \end{pmatrix} \in \mathcal D,
\]
 and the channels $C$ and $\bar C$ as before. This time the maximal expected utilities are
\[
\max_{D \in \mathit{\Phi}(C)}\tr {U\cdot D \cdot \mathit{\Pi}}= \nicefrac{1}{2}, \quad 
\max_{D \in \mathit{\Phi}\left(\bar C\right)}\tr {U\cdot D \cdot \mathit{\Pi}}=\nicefrac{2}{3},
\]
respectively. Thus, $C$ is not more Blackwell-useful with respect to $\mathcal D$ than $\bar C$. \qed
\end{itemize}
\end{example}

Since it seems hard to extend Proposition \ref{prop:exact}, we now focus on the necessity of the Shannon-order for the reduced Blackwell-usefulness.
Let us examine for which subsets ${\mathcal U}$ of utility matrices the reverse implication holds:
\[
C \trianglerighteq_S \bar C \quad 
\Leftarrow \quad
C \succcurlyeq_B^{{\mathcal U}} \bar C.
\]
It follows due to Remark \ref{rem:indif} that the set $\mathcal I$ is not necessary for that. Otherwise, all channels would be Shannon-garblings of each other, trivially a false statement. It turns out that the necessity also fails for the set $\mathcal E$. 

\begin{example}[Failure of necessity of $\mathcal E$]
\label{ex:exactsmall}
The reverse implication in (\ref{eq:buso}) does not  hold for ${\mathcal U}=\mathcal E$ in case of $m=n=2$.
%The set of exact utility matrices can be written as
%\[
%\mbox{$\mathcal E$}=\left\{ U \left\vert \mbox{There exists } a>0 \mbox{ with } U=a\cdot P\begin{pmatrix}1&0\\0&1\end{pmatrix} \mbox{ or } U=a\cdot\begin{pmatrix}0&1\\1&0\end{pmatrix} \right. \right\}.
%\]
We consider two channels
\[
C=\begin{pmatrix}1 & \nicefrac{1}{2} \\ 0 & \nicefrac{1}{2}\end{pmatrix}, \quad
\bar C=\begin{pmatrix}\nicefrac{1}{4} & \nicefrac{3}{4} \\  \nicefrac{3}{4} & \nicefrac{1}{4}\end{pmatrix}.
\] 
Neither of these channels is a Shannon-garbling of the other, as it can be seen from a straight-forward calculation. Nevertheless, $C$ is more Blackwell-useful with respect to $\mathcal E$ than $\bar C$. In fact, for all $U = a \cdot P \in \mathcal E$ with an arbitrary, but fixed permutation matrix $P$ it holds:
\[
\max_{D \in \mathit{\Phi}(C)}\tr{U\cdot D \cdot \mathit{\Pi}}
=\max_{D \in \mathit{\Phi}\left(\bar C \right)}\tr{U\cdot D \cdot \mathit{\Pi}}=a \cdot \nicefrac{3}{4}.
\]\qed
\end{example}

Hence, neither $\mathcal E$ nor its subsets can provide the equivalence between Shannon-order and reduced Blackwell-usefulness.
As it turns out, the equivalence (\ref{eq:buso}) does not hold in general for any subset $U$ of utility matrices. This is shown in Theorem \ref{th:buso} for the case $n=2^{m-2}+1$, $m \geq 2$. Lemma \ref{lem:2x2} treats first the case $n=m=2$.

\begin{lemma}[$n=2, m=2$]
\label{lem:2x2}
 The Shannon-order and reduced Blackwell-usefulness are not equivalent for any subset $U \subseteq \R^{2 \times 2}$ of utility matrices.
\end{lemma}

\proof
Let us suppose that there exists a subset ${\mathcal U}$, such that for every two channels $C, \bar C \in \R^{2 \times 2}$ with the same input alphabet it holds:
\[
C \trianglerighteq_S \bar C \quad 
\Leftrightarrow \quad
C \succcurlyeq_B^{\mathcal U} \bar C.
\]
 ${\mathcal U}$ is a subset of $2\times 2$-matrices which can be written as
\[
\begin{pmatrix} a & a + \varepsilon_1 \\
b + \varepsilon_2 & b\end{pmatrix}
\]
for some $a,b,\varepsilon_1,\varepsilon_2 \in \R$.
Due to Remark \ref{rem:indif}, the addition of indifferent utility matrices does not affect the above equivalence. Hence, we can assume without loss of generality that all matrices in $\mathcal U$ can be written as
\[ U= \begin{pmatrix} a & a + \varepsilon_1 \\
b + \varepsilon_2 & b\end{pmatrix} + \begin{pmatrix} -a & -a \\ -b & -b\end{pmatrix} = \begin{pmatrix}0 & \varepsilon_1 \\ \varepsilon_2 & 0\end{pmatrix}.
\]
Each of those matrices belongs to one of the following sets:
\[
\begin{array}{rcl}
    {\mathcal U}_{\le} &=&\left\{ U \in {\mathcal U} \, \left\vert \,
U=\begin{pmatrix}0 & \varepsilon_1 \\ \varepsilon_2 & 0\end{pmatrix}, \varepsilon_1 \cdot \varepsilon_2 \le 0 \right.\right\}, \\ \\
{\mathcal U}_{\ne} &=&\left\{ U \in {\mathcal U} \,\left\vert \,
U=\begin{pmatrix}0 & \varepsilon_1 \\ \varepsilon_2 & 0\end{pmatrix}, \varepsilon_1 \cdot \varepsilon_2 > 0, \varepsilon_1 \ne \varepsilon_2 \right.\right\},\\ \\
{\mathcal U}_{=}&=&\left\{ U \in {\mathcal U} \,\left\vert \,
U=\begin{pmatrix}0 & \varepsilon_1 \\ \varepsilon_2 & 0\end{pmatrix}, \varepsilon_1 \cdot \varepsilon_2 > 0, \varepsilon_1 = \varepsilon_2 \right.\right\}.
\end{array}
\]
The subset $\mathcal U_{\le}$ only contains matrices with a dominant column. Let a decision maker have a utility matrix $U \in \mathcal U_{\le}$ with the $k$-th dominant column. Then, maximal expected utility is gained by taking the action $A= \left(a_{ij}\right)$ with $a_{kj} = 1$ for all $j$, and $a_{ij}=0$ otherwise.
Since the maximal expected utility is the same for any channel, 
we can assume without loss of generality that the subset ${\mathcal U}_{\le}$ is empty.

Now, consider a utility matrix $U \in {\mathcal U}_{\ne}$ and the two channels
\[
C_1=\begin{pmatrix}1 & \nicefrac{1}{2} \\ 0 & \nicefrac{1}{2} \end{pmatrix}  , \qquad
C_2=\begin{pmatrix}\nicefrac{1}{2} & 1 \\ \nicefrac{1}{2} & 0 \end{pmatrix}.
\]
Each of those channels is a Shannon-garbling of the other with the matrices
\[
M=\begin{pmatrix}1 & 0\\ 0 & 1\end{pmatrix} , \qquad N=\begin{pmatrix} 0 & 1 \\ 1 & 0 \end{pmatrix}.
\]
%Hence, for every $U\in \mathcal U$ the maximal expected utility has to be the same for both channels.
But, for $\varepsilon_1 > \varepsilon_2 > 0$ we have
\[
\begin{array}{rcl}
\displaystyle \max_{D \in \mathit{\Phi}(C_1)}\tr{U\cdot D \cdot \mathit{\Pi}} &=& 
\displaystyle \nicefrac{\displaystyle \varepsilon_1}{2} + \nicefrac{\displaystyle \varepsilon_2}{4}, \\ \\ 
\displaystyle \max_{D \in \mathit{\Phi}\left(C_2 \right)}\tr{U\cdot D \cdot \mathit{\Pi}} &=& \displaystyle \nicefrac{\displaystyle \varepsilon_1}{2} + \max\{0, \nicefrac{\displaystyle \varepsilon_2}{2}-\nicefrac{\displaystyle \varepsilon_1}{4}\}.
\end{array}
\]
For $\varepsilon_2 > \varepsilon_1 > 0$ we have
\[
\begin{array}{rcl}
\displaystyle \max_{D \in \mathit{\Phi}(C_1)}\tr{U\cdot D \cdot \mathit{\Pi}} &=& \displaystyle
\nicefrac{\displaystyle \varepsilon_1}{2} + \nicefrac{\displaystyle \varepsilon_2}{4} + \max\{0, \nicefrac{\displaystyle \varepsilon_2}{4}-\nicefrac{\displaystyle \varepsilon_1}{2}\}, \\ \\
\displaystyle 
\max_{D \in \mathit{\Phi}\left(C_2 \right)}\tr{U\cdot D \cdot \mathit{\Pi}} &=& \displaystyle \nicefrac{\displaystyle \varepsilon_2}{2} + \nicefrac{\displaystyle \varepsilon_1}{4}.
\end{array}
\]
The remaining cases lead analogously to 
\[
\max_{D \in \mathit{\Phi}(C_1)}\tr{U\cdot D \cdot \mathit{\Pi}} \not =
\max_{D \in \mathit{\Phi}\left(C_2 \right)}\tr{U\cdot D \cdot \mathit{\Pi}}.
\]
Therefore, ${\mathcal U}_{\ne}$ is to be empty.

Altogether, we have $\mathcal U = {\mathcal U}_{=}$, and every $U \in \mathcal{U}$ can be written as
\[
U=\begin{pmatrix}0&\varepsilon\\ \varepsilon&0\end{pmatrix}
\]
with some $\varepsilon\ne 0$.
By using Remark \ref{rem:indif} again, we may add an indifferent matrix to conclude that every $U \in \mathcal{U}$ can be written as
\[
U=\begin{pmatrix}\varepsilon&0 \\ 0 & \varepsilon\end{pmatrix}
\qquad \mbox{or} \qquad
U=\begin{pmatrix} 0 & \varepsilon \\ \varepsilon& 0\end{pmatrix},
\]
where $\varepsilon>0$.
Thus, $\mathcal U$ is a subset of $\mathcal E$. However, due to Example \ref{ex:exactsmall}, the subsets of $\mathcal E$ do not provide the equivalence (\ref{eq:buso}).
\qed

Lemma \ref{lem:2x2} can be generalized as follows.

\begin{theorem}[$n=2^{m-2}+1, m \geq 2$]
\label{th:buso}
 The Shannon-order and reduced Blackwell-usefulness are not equivalent for any subset $U \subseteq \R^{n \times m}$ of utility matrices with $n=2^{m-2}+1$, $m \geq 2$.
\end{theorem}

\proof
%%%%%%%%%%%%%%%%%%%%%%%%%
For $m=2$ this is due to Lemma \ref{lem:2x2}.
Let us suppose that there exists a subset ${\mathcal U}\subseteq \R^{n\times m}, m\ge 3$, such that for every two channels $C, \bar C \in \R^{m \times n}$ with the same input alphabet it holds:
\[
C \trianglerighteq_S \bar C \quad 
\Leftrightarrow \quad
C \succcurlyeq_B^{\mathcal U} \bar C.
\]
%First, we define $k=2^{m-3}+1$.
%, \qquad l=n-k=2^{m-3}.\]
Let $U \in \mathcal U$ be a fixed utility matrix. For all permutation matrices $P \in \R^{n \times n}$ and all channels $C \in \R^{m \times n}$ it holds:
\[
\max_{D \in \mathit{\Phi}(C)}\tr{U\cdot D \cdot \mathit{\Pi}}=\max_{D \in \mathit{\Phi}(C)}\tr{(P\cdot U)\cdot D \cdot \mathit{\Pi}}
\]
This is due to  the fact, that for every channel $C \in \R^{m\times n}$ it holds both:
\[
C \trianglerighteq_S C \cdot P\qquad \mbox{and}\qquad C\cdot P \trianglerighteq_S  C.
\]
Therefore, for every $U\in \mathcal{U}$ the subset $\left\{P \cdot U \in \mathcal{U}\, \left\vert\, P \mbox{ is a permutation matrix}\right. \right\}$  can be replaced by a single utility matrix $P\cdot U$.
Next, we want to chose an appropriate permutation matrix $P$ in dependence on $U$.
Since every utility matrix $U=\left(u_{ij}\right)\in \mathcal{U}$ has $2^{m-2}+1$ rows, there exist at least $k=2^{m-3}+1$ indices $i_1,\ldots,i_{k}$ such that it holds: 
\[
u_{i1}\le u_{i2} \quad \mbox{for } i=i_1,\ldots,i_{k}
\quad
\mbox{or}
\quad
u_{i1}\ge u_{i2} \quad \mbox{for } i=i_1,\ldots,i_{k}.
\]
%For every $U$ it exists a permutation matrix $P$ such that for $P\cdot U$ the corresponding indices are $i_1=1, \ldots, i_k=k $.
If $i_1,\ldots,i_{k}$ can be chosen as $1,\ldots, k$, we say $U$ fulfills the ordering condition.
For every $U$ it exists a permutation $P$ such that for $P\cdot U$ the latter holds.
Hence, we assume without loss of generality that $\mathcal{U}$ consists only of utility matrices which fulfill the ordering condition.

We now define a new set of utility matrices
$\mathcal{\bar U}\subseteq \R^{k \times (m-1)}$ which consists of  
\[
\bar U= \begin{pmatrix} \max \left\{u_{11}, u_{12}\right\}&u_{13}& \ldots  & \ u_{1m}\\
\vdots &\vdots& & \vdots\\
\max \left\{ u_{k1}, u_{k2}\right\}&u_{k3}& \ldots & u_{km}\end{pmatrix}
\]
for some utility matrix $U = \left(u_{ij}\right)\in\mathcal{U}$. We consider two channels 
$Z=\left(z_{ij}\right), \bar Z=\left(\bar z_{ij}\right)\in \R^{(m-1)\times k}$. Let us first assume that $\bar Z$ is a Shannon-garbling of $Z$, i.\,e.
there exist stochastic matrices $M, N$ with
\[
\bar Z= M \cdot Z \cdot N.
\]
We define the channels $C, \bar C \in \R^{m \times n}$:
\[
C:=\begin{pmatrix}
z_{11}&\ldots&z_{1k}&0&\ldots&0\\
\vdots&&\vdots&\vdots&&\vdots\\
 z_{(m-1)1}&\ldots&z_{(m-1)k}&0&\ldots&0\\
0& \ldots &0&1&\ldots&1\end{pmatrix},
\]
\[
\bar C:=\begin{pmatrix}
\bar z_{11}&\ldots&\bar z_{1k}&0&\ldots&0\\
\vdots&&\vdots&\vdots&&\vdots\\
 \bar z_{(m-1)1}&\ldots&\bar z_{(m-1)k}&0&\ldots&0\\
0& \ldots &0&1&\ldots&1\end{pmatrix}.
\]
Then it holds:
\[
\bar C=\begin{pmatrix}
M&\mathbf{0}\\
\mathbf{0}&1\end{pmatrix}\cdot C\cdot \begin{pmatrix}N&\mathbf{0}\\ \mathbf{0}&\mathbf{I}\end{pmatrix},
\]
where $M$ and $N$ are the above mentioned matrices, $\mathbf{I}$ represents the identity matrix, and  $\mathbf{0}$ denotes  zero matrices of proper dimensions.
Therefore, it holds $C \trianglerighteq_S \bar C$, hence, also $C \succcurlyeq_B^{\mathcal U} \bar C$. 
Next, we rewrite the maximal expected utility of $C=\left(c_{ij}\right)$ when endowed with $U=\left(u_{ij}\right) \in \mathcal{U}$:
\[
\max_{D \in \mathit{\Phi}(C)}\tr{U\cdot D \cdot \mathit{\Pi}}
=  \frac{1}{n}\cdot \max_{A \in \mathcal{A}}\tr{C\cdot U\cdot A}\notag. \\
%&=  \frac{1}{n}\cdot \max_{A \in \mathcal{A}}\tr{\begin{pmatrix}z_{11}\cdot u_{11}\ldots&z_{11}\cdot u_{1k}&0&\ldots&0\\
%\vdots&&\vdots&\vdots&&\vdots\\
%z_{(m-1)1}\cdot u_{11}\ldots&z_{(m-1)1}\cdot u_{1k}&0&\ldots&0\end{pmatrix}\cdot A}
\]
We decompose $\tr{C\cdot U \cdot A}$ into the following sum:
\[
\begin{array}{rcl}
\tr{C\cdot U \cdot A}&=&\left(z_{11}\cdot u_{11}+\ldots+z_{1k}\cdot u_{k1}\right)\cdot a_{11}\\  \\
&&+\left(z_{11}\cdot u_{12}+\ldots+z_{1k}\cdot u_{k2}\right)\cdot a_{21}\\  \\
&&+\ldots\\ \\
&&+\left(z_{11}\cdot u_{1m}+\ldots+z_{1k}\cdot u_{km}\right)\cdot a_{m1}\\  \\
&&+\ldots \\ \\
&&+ \left(z_{(m-1)1}\cdot u_{11}+\ldots+z_{(m-1)k}\cdot u_{k1}\right)\cdot a_{1(m-1)}\\  \\
&&+ \left(z_{(m-1)1}\cdot u_{12}+\ldots+z_{(m-1)k}\cdot u_{k2}\right)\cdot a_{2(m-1)}\\  \\
&&+ \ldots \\ \\
&&+ \left(z_{(m-1)1}\cdot u_{1m}+\ldots+z_{(m-1)k}\cdot u_{km}\right)\cdot a_{m(m-1)}\\  \\
&&+ \sum\limits_{i=k+1}^n u_{i1} \cdot a_{1m}+\ldots+ \sum\limits_{i=k+1}^n u_{im} \cdot a_{mm}.\\ \\
\end{array}
\]
Since $U$ fulfills the ordering criteria we can now simplify this. Without loss of generality we assume that it holds:
\[
u_{i1}\le u_{i2}\quad \mbox{for }i=1,\ldots,k.
\]
The remaining case can be namely proven analogously.
Hence, we obtain:
\[
\begin{array}{rcl}
\tr{C\cdot U \cdot A}&\le&
\left(z_{11}\cdot u_{12}+\ldots+z_{1k}\cdot u_{k2}\right)\cdot \left(a_{11} + a_{21}\right)\\  \\
&&+\ldots\\ \\
&&+\left(z_{11}\cdot u_{1m}+\ldots+z_{1k}\cdot u_{km}\right)\cdot a_{m1}\\  \\
&&+\ldots \\ \\
&&+ \left(z_{(m-1)1}\cdot u_{12}+\ldots+z_{(m-1)k}\cdot u_{k2}\right)\cdot \left(a_{1(m-1)} + a_{2(m-1)}\right)\\  \\
&&+ \ldots \\ \\
&&+ \left(z_{(m-1)1}\cdot u_{1m}+\ldots+z_{(m-1)k}\cdot u_{km}\right)\cdot a_{m(m-1)}\\  \\
&&+ \displaystyle \max_{j=1,\ldots,m}\sum\limits_{i=k+1}^n u_{ij}.
\end{array}
\]
We define $B=(b_{ij})\in \R^{(m-1) \times (m-1)}$ by setting
\[
  \begin{array}{lcll}
  b_{1j} &=& a_{1j} + a_{2j} & \mbox{for } j =1, \ldots, m-1, \\
  b_{ij} &=& a_{i+1,j} &\mbox{for } i=2, \ldots, m-1, j =1, \ldots, m-1.
  \end{array}
\]
We also set
\[
  s = \max_{j=1,\ldots,m}\sum\limits_{i=k+1}^n u_{ij}.
\]
Then, the above inequality becomes
\[
  \tr{C\cdot U \cdot A} \leq \tr{Z \cdot \bar U \cdot B} +s.
\]
In particular, it follows that
\[
  \max_{A \in \mathcal{A}} \tr{C\cdot U \cdot A} = \max_{B \in \mathcal{A}} \tr{Z \cdot \bar U \cdot B} +s.
\]
Analogously:
\[
  \max_{A \in \mathcal{A}} \tr{\bar C\cdot U \cdot A} = \max_{B \in \mathcal{A}} \tr{\bar Z \cdot \bar U \cdot B} +s.
\]
Overall, we have proved that
\[
  Z \trianglerighteq_S \bar Z \quad 
\Rightarrow \quad
Z \succcurlyeq_B^{\bar{\mathcal{U}}} \bar Z.
\]

Now, let us assume that $Z \succcurlyeq_B^{\bar{\mathcal{U}}} \bar Z$. Then, using the construction above, we have $C \succcurlyeq_B^{\mathcal{U}} \bar C$, hence, $C \trianglerighteq_S \bar C$. With stochastic matrices $\bar M, \bar N$ we have:
\[
\bar C= \bar M \cdot C \cdot \bar N.
\]
Let us write by using blocks of appropriate size:
\[
  \bar M = \left(\begin{array}{cc}
       \bar M_{11} &  \bar M_{12} \\
      \bar M_{21} & \bar M_{22}
  \end{array} \right), \quad 
  \bar N = \left(\begin{array}{cc}
       \bar N_{11} & \bar N_{12} \\
       \bar N_{21} & \bar N_{22}
  \end{array}\right).
\]
These matrices are of the following dimensions
\[
 \begin{array}{llll}
\bar M_{11} \in \R^{(m-1)\times(m-1)}, &
\bar M_{12} \in \R^{m-1}, &
\bar M_{21}^T \in \R^{m-1}, &
\bar M_{22} \in \R, \\ \\
\bar N_{11} \in \R^{k\times k}, &
\bar N_{12} \in \R^{k\times(n-k)}, &
\bar N_{21} \in \R^{(n-k)\times k}, &
\bar N_{22} \in \R^{(n-k)\times(n-k)}.
 \end{array}
\]
By multiplying out, we obtain:
\[
\begin{array}{rcl}
   \bar Z &=& \bar M_{11} \cdot Z \cdot\bar N_{11} + \bar M_{12} \cdot E \cdot \bar N_{21}, \\ \\
   \mathbf{0} &=& \bar M_{21} \cdot Z \cdot\bar N_{11} + \bar M_{22} \cdot E \cdot\bar N_{21}, \\ \\
   \mathbf{0} &=& \bar M_{11} \cdot Z \cdot\bar N_{12} + \bar M_{12} \cdot E \cdot\bar N_{22}, \\ \\
   E &=& \bar M_{21} \cdot Z \cdot\bar N_{12} + \bar M_{22} \cdot E \cdot\bar N_{22},     
\end{array}
\]
where by $E^T\in \R^{n-k}$ we denote the vector of ones.

\underline{\textsc{case 1:}}\;$\bar M_{12}= \mathbf{0}$.
\newline
Then, by using the first and the second equation, we have:
\[
   \bar Z =  \left( \bar M_{11} + \left(\begin{array}{c}
        \mathbf{0}  \\
        \bar M_{21}
   \end{array}\right)\right) \cdot Z \cdot \bar N_{11},
\]
where the matrix $\left( \bar M_{11} + \left(\begin{array}{c}
        \mathbf{0}  \\
        \bar M_{21}
   \end{array}\right)\right)$ is stochastic. From the second equation we have $\bar N_{21}=\mathbf{0}$, hence $\bar N_{11}$ is also stochastic.
  
\underline{\textsc{case 2:}}\;$\bar N_{21}= \mathbf{0}$.
\newline
Then, by using the first and the second equation, we have:
\[
   \bar Z =  \left( \bar M_{11} + \left(\begin{array}{c}
        \mathbf{0}  \\
        \bar M_{21}
   \end{array}\right)\right) \cdot Z \cdot \bar N_{11},
\]
where the matrices $\left( \bar M_{11} + \left(\begin{array}{c}
        \mathbf{0}  \\
        \bar M_{21}
   \end{array}\right)\right)$ and $\bar N_{11}$ are stochastic.

\underline{\textsc{case 3:}}\;$\bar M_{12} \not = \mathbf{0}, \bar N_{21} \not = \mathbf{0}$.
\newline
From the second equation follows that $\bar M_{22} = \mathbf{0}$,
and from the third equation that $\bar N_{22} = \mathbf{0}$.
Hence, from the fourth equation follows that $\bar M_{21} \cdot Z \cdot\bar N_{12}=E$. Thus,  the first equation becomes:
\[
\begin{array}{rcl}
\bar Z &=& \bar M_{11} \cdot Z \cdot\bar N_{11} + \bar M_{12} \cdot \bar M_{21} \cdot Z \cdot\bar N_{12} \cdot \bar N_{21}\\ \\
&=& \bar M_{11} \cdot Z \cdot\bar N_{11} + \bar M_{12} \cdot \bar M_{21} \cdot Z \cdot\bar N_{12} \cdot \bar N_{21} \\ \\ &&+  \bar M_{12} \cdot \underbrace{\bar M_{21} \cdot Z \cdot\bar N_{11}}_{ \begin{array}{c}
     =\mathbf{0} \\ \mbox{due to the} \\ \mbox{second equation}
\end{array}}+ \underbrace{\bar M_{11} \cdot Z \cdot\bar N_{12}}_{ \begin{array}{c}
     =\mathbf{0} \\ \mbox{due to the} \\ \mbox{third equation}
\end{array}} \cdot \bar N_{21}\\ \\
&=& \left( \bar M_{11}+  \bar M_{12} \cdot \bar M_{21} \right) \cdot Z \cdot \left(\bar N_{11} + \bar N_{12} \cdot \bar N_{21}\right).
\end{array} 
\]
Due to stochasticity of $\bar M_{12}$ and $ \bar N_{12}$, it follows that 
$\left( \bar M_{11}+  \bar M_{12} \cdot \bar M_{21} \right)$ and $\left(\bar N_{11} + \bar N_{12} \cdot \bar N_{21}\right)$ are stochastic.

From these cases, we conclude that $\bar Z$ is a Shannon-garbling of $Z$, i.\,e.
\[
   Z \trianglerighteq_S \bar Z \quad 
\Leftarrow \quad
Z \succcurlyeq_B^{\bar{\mathcal{U}}} \bar Z.
\]

Overall, we reduced the dimension from $n=2^{m-2}+1$ to $k=2^{(m-1)-2}+1$. The contradiction follows by induction and Lemma \ref{lem:2x2}.
\qed
%%%%%%%%%%%%%%%%%%%%%%%%

%The failure of equivalence of Shannon-order and reduced Blackwell-usefulness can be generalized by using the following Lemma \ref{lem:general}.

%\begin{lemma}
%\label{lem:general}
%If the Shannon-order and reduced Blackwell-usefulness are equivalent for a subset $U \subseteq \R^{n \times m}$ of utility matrices, then it is also true 
%for a subset $\bar U \subseteq \R^{(n-1) \times m}$.
%\end{lemma}
%\proof
%See Appendix \ref{app}. \qed

%\begin{theorem}
%\label{th:buso}
%The Shannon-order and reduced Blackwell-usefulness are not equivalent for any subset $U \subseteq \R^{n \times m}$ of utility matrices with $m=2,3$ and any $n\geq m$.
%\end{theorem}

%\proof

Theorem \ref{th:buso} states that there does not exist a set of utility matrices, such that the corresponding reduced Blackwell-usefulness and Shannon-order are equivalent, at least for $n=2^{m-2}+1$ and $m \geq 2$. Note that for every nontrivial channel $C$ we may assume that $C \in \R^{m\times (2^{m-2}+1)}$ with some $m \geq 2$. This is achievable by duplicating input letters, while adjusting the input distribution accordingly or adding output letters, which will not be reported at all.  

Recall that Theorem \ref{th:buso} holds for channels with uniform distribution of the input alphabet. We will now show that it also holds for arbitrary distributions of the input alphabet.

\begin{remark}[General distribution of input alphabet]
\label{re:allpi}
Let the input alphabet be generally distributed with probabilities $\bar \pi_i > 0$, $i=1, \ldots, n$, i.\,e. 
\[
\mathit{\bar \Pi}=\mbox{diag}\left(\bar \pi_1, \ldots, \bar \pi_n\right).
\]
We claim that there does not exist a subset $\bar{\mathcal U}$ of utility matrices such that
\[
C \trianglerighteq_S \bar C \quad 
\Leftrightarrow \quad
C \succcurlyeq_B^{\bar{\mathcal U}} \bar C
\]
holds for every two channels $C, \bar C \in \R^{m \times n}$ with the same input alphabet.
In fact, let us assume on the contrary that such a subset $\bar{\mathcal U}$ of utility matrices exists. 
We again denote by $\mathit{\Pi}$ the uniform distribution, and define:
\[
{\mathcal U}:=\left\{ U \, \left\vert \, \mbox{There exists } \bar U \in \bar{\mathcal U} \mbox{ with } U= \mathit{ \Pi}^{-1}\cdot \mathit{ \bar \Pi} \cdot \bar U \right. \right\}.
\]
Hence, 
\[
\max_{D \in \mathit{\Phi}(C)}\tr {\bar U\cdot D \cdot \mathit{\bar \Pi}}\ge \max_{D \in \mathit{\Phi}\left(\bar C\right)}\tr {\bar U\cdot D \cdot \mathit{\bar \Pi}}
\]
holds for all $\bar U \in \bar{\mathcal U}$ if and only if
\[
\max_{D \in \mathit{\Phi}(C)}\tr { U\cdot D \cdot \mathit{\Pi}}\ge \max_{D \in \mathit{\Phi}\left(\bar C\right)}\tr {U\cdot D \cdot \mathit{\Pi}}
\]
holds for all $ U \in {\mathcal U}$.
Thus, for channels $C, \bar C \in \R^{m \times n}$ with uniformly distributed input alphabets we have:
\[
C \trianglerighteq_S \bar C \quad 
\Leftrightarrow \quad
C \succcurlyeq_B^{{\mathcal U}} \bar C.
\]
From Theorem \ref{th:buso} we know that the latter equivalence is not valid for any subset ${\mathcal U}$ of utility matrices at least for $n=2^{m-2}+1$ and $m\geq 2$, a contradiction. \qed
\end{remark}

\section{Convexified Shannon-usefulness}

We have seen in Section \ref{sec:neg} that the Shannon-order is not equivalent to (reduced) Blackwell-usefulness. To overcome this difficulty we instead characterize the convexified Shannon-order by an appropriate notion of usefulness.

\begin{definition}[Convexified Shannon-order,  \citet{shannon:1958}]
Let  $C, \bar C \in \R^{m \times n}$ be two channels with the same input alphabet. We say $\bar C$ is a convexified Shannon-garbling of $C$ (denoted by $C \trianglerighteq_{cS} \bar C$) if there exist a probability distribution
$q_j$, $j=1,\ldots,\ell$, and stochastic matrices
$M_j \in \R^{m \times m}$, $N_j \in \R^{n \times n}$, $j=1,\ldots,\ell$,
with
\[
\bar C= \sum\limits_{j=1}^\ell q_j\cdot M_j \cdot C \cdot N_j.
\]
We call $\trianglerighteq_{cS}$ the partial convexified Shannon-order of channels.
\end{definition}

Now we introduce the corresponding notion of the convexified Shannon-usefulness.
For that, let the convexified Shannon policy space be defined as follows:
\[
\mathit{\Phi}_{cS}(C) = \left\{D \in \R^{m \times n} \,\left\vert \,
\begin{array}{l}
\mbox{There exist a probability distribution }
p_i, i=1,\ldots,k, \\ \mbox{and } A_i, B_i \in {\mathcal A}, i=1, \ldots, k
\mbox{ with } D=\sum\limits_{i=1}^k p_i \cdot A_i \cdot C \cdot B_i 
\end{array}  \right.
\right\}.
\]
Note that the convexified Shannon policy space $\mathit{\Phi}_{cS}(C)$ consists of all convexified Shannon-garblings of $C$.

\begin{lemma} 
\label{lem:cc}
The convexified Shannon policy space $\mathit{\Phi}_{cS}(C)$ is convex and compact.
\end{lemma}
\proof
From definition it immediately follows that $\mathit{\Phi}_{cS}(C)$ is convex.
In order to prove that $\mathit{\Phi}_{cS}(C)$ is compact, we show that it is the convex hull of a finite set. For that, let $A\in \R^{m\times m}$ be a stochastic matrix and $\mathcal{L} \subset \R^{m\times m}$ the set of matrices whose columns are coordinate vectors. Since this set is finite, we can write it as $\mathcal{L}=\left\{L_1,\ldots L_{m^m}\right\}$. Thus, there exist $\alpha_j \ge 0$, $j=1,\ldots,m^m$ such that it holds:
\[
A= \sum\limits_{j=1}^{m^m} \alpha_j \cdot L_j, \qquad \sum\limits_{j=1}^{m^m} \alpha_j=1.
\]
Analogously, let $B\in \R^{n\times n}$ be a stochastic matrix and $\mathcal{R} \subset \R^{n\times n}$ the set of matrices whose columns are coordinate vectors. 
We set  $\mathcal{R}=\left\{R_1,\ldots R_{n^n}\right\}$. Thus, there exist $\beta_{\ell} \ge 0$, $\ell=1,\ldots,n^n$ such that it holds:
\[
B= \sum\limits_{\ell=1}^{n^n} \beta_{\ell} \cdot R_{\ell}, \qquad \sum\limits_{\ell=1}^{n^n} \beta_{\ell}=1.
\]
Hence, we have:
\[
A \cdot C \cdot B=  \sum\limits_{j=1}^{m^m} \sum\limits_{\ell=1}^{n^n}\alpha_j \cdot \beta_{\ell} \cdot L_j \cdot C \cdot  R_\ell.
\]
Therefore, every element $D \in \mathit{\Phi}_{cS}(C)$ can be written as
\[
\begin{array}{rcl}
D= \sum\limits_{i=1}^k p_i \cdot A_i \cdot C \cdot B_i  &=&\sum\limits_{i=1}^k p_i \cdot  \sum\limits_{j=1}^{m^m} \sum\limits_{\ell=1}^{n^n}\alpha_{ij} \cdot \beta_{i\ell} \cdot L_j \cdot C \cdot  R_\ell\\ \\
&=&\sum\limits_{j=1}^{m^m} \sum\limits_{\ell=1}^{n^n}\gamma_{j\ell} \cdot \left( L_j \cdot C \cdot  R_\ell\right)
\end{array}
\]with 
\[
\gamma_{j\ell}=\sum\limits_{i=1}^k p_i \cdot \alpha_{ij} \cdot \beta_{i\ell},\qquad \sum\limits_{j=1}^{m^m} \sum\limits_{\ell=1}^{n^n}\gamma_{j\ell}=1.
\]
In particular, it follows that
\[
D \in \mbox{Conv}\left(\left\{ L \cdot C \cdot  R \; \left\vert \, L \in \mathcal{L}, R \in \mathcal{R}\right. \right\} \right).
\]
Thus, the convexified Shannon policy space is a subset of the latter convex hull, i.\,e. 
\[
\mathit{\Phi}_{cS}(C) \subseteq  \mbox{Conv}\left(\left\{ L \cdot C \cdot  R \; \left\vert \, L \in \mathcal{L}, R \in \mathcal{R}\right. \right\} \right).
\]
It is easy to see that the reverse also holds. In fact, $L \cdot C \cdot  R \in \mathit{\Phi}_{cS}(C)$ for in particular stochastic matrices $L \in \mathcal{L}$, $R \in \mathcal{R}$. The convexity of $\mathit{\Phi}_{cS}(C)$ provides the assertion.
\qed
%First we want to introduce, the term of Shannon-usefulness. Since Blackwell-usefulness was a optimization-problem over the policy space $\mathit{\Phi}(C)$ of a channel $C$. We also want to introduce the Shannon policy space.
%Furthermore, we want to convexify both definitions. This will be useful in order to proof the final theorem for the Shannon-order. Convexifying leads also to a more general term of channel inclusion than the Shannon-order, which we will call the convexified Shannon-order \citep{shannon:1958}.

\begin{definition}[Convexified Shannon-usefulness]
Let $C, \bar C \in \R^{m\times n}$ be two channels with the same input-alphabet.
We say that $C$ is more convexified-Shannon-useful than $\bar C$ (denoted by $C \succcurlyeq_{cS} \bar C$) if for all utility matrices $U \in \mathbb{R}^{n \times m}$ it holds:
\[
\max_{D \in \mathit{\Phi}_{cS}(C)}\tr {U\cdot D \cdot  \mathit{ \Pi}}\ge \max_{D \in \mathit{\Phi}_{cS}\left(\bar C\right)}\tr {U\cdot D \cdot  \mathit{ \Pi}}.
\]
This means that every decision maker gains by using $C$ at least the utility he  or she would gain by $\bar C$.
\end{definition}

\begin{remark}[Interpretation of convexified Shannon-usefulness]
When evaluating a channel by Blackwell-usefulness, decision makers optimize over their reactions to the channel's output. However, Shannon-usefulness endows decision makers with more possibilities. First, decision makers are not only allowed to react to the channel's output, but also to code the channel's input. This corresponds to multiplication of the channel from the right by a stochastic matrix. Thus, the decision makers are able to distribute the noise of the channel to any input letter. The columns of the original channel can be interpreted as the output distribution of a given input letter. The decision makers' choice is to assign to an every input letter a desirable output distribution. They are allowed to collate different output distributions by forming their convex combinations. Moreover, they may replace the output distributions of some input letters by those of the others. This means, before using the actual channel, decision makers are allowed to code the original message. After this, they will react to the output of the coded channel.
This corresponds to multiplication of the channel from the left by a stochastic matrix.
Additionally, decision makers are allowed to repeat this process by using various coding protocols and determining other reactions accordingly. The repetition is due to a probability distribution. This corresponds to convexification. Finally, they will optimize over all coding protocols, possible reactions and probability distributions.
Practically this means that, when decision makers react to a received message over a noisy channel, the convexified Shannon-usefulness enables agreements on coding, as well as on sending the message repeatedly.
\qed

\end{remark}

In a similar way to the proofs of Blackwell's Theorem by \cite{leshno:1992}, and \cite{perez-richet:2017}, we show that the convexified Shannon-usefulness characterizes the convexified Shannon-order.

\begin{theorem}
\label{th:shannonstheorem}
It  holds for channels $C, \bar C \in \R^{m \times n}$ with the same input alphabet:
\[
C \trianglerighteq_{cS} \bar C \quad 
\Leftrightarrow \quad
C \succcurlyeq_{cS} \bar C.
\]
\end{theorem}
\proof
%Let $\Pi$ be the distribution matrix of the input alphabet.

\underline{\textsc{step 1:}} \;$C\succcurlyeq_{cS} \bar C \Rightarrow \mathit{\Phi}_{cS}(C) \supseteq \mathit{\Phi}_{cS}(\bar C).$\\
%First, we will show that it follows from $C \succcurlyeq_{cS} \bar C$, that the convexified Shannon policy space of $\bar C$ is a subset of that of $C$, i.e. $\mathit{\Phi}_{cS}\left(\bar C\right)\subseteq \mathit{\Phi}_{cS}\left(C\right)$.
Suppose on the contrary there exists $\bar D = \left(\bar d_{ij}\right)\in \mathit{\Phi}_{cS}\left(\bar C\right) \backslash \mathit{\Phi}_{cS}\left(C\right)$.
Since $\mathit{\Phi}_{cS}\left(C\right)$ is closed and convex due to Lemma \ref{lem:cc}, we may apply the separation theorem. Hence, there exists a linear functional $U=\left(u_{ij}\right) \in \mathbb{R}^{m \times n}$ such that for all
$D =\left(d_{ij}\right)\in \mathit{\Phi}_{cS}\left( C\right)$ it holds:
\[
\sum\limits_{i=1}^{m}\sum\limits_{j=1}^{n}{ u_{ij} \cdot \bar d_{ij}} > \sum\limits_{i=1}^{m}\sum\limits_{j=1}^{n}{u_{ij} \cdot d_{ij}}.
\]
We define:
\[
\bar U:= \Pi^{-1}\cdot U^T.
\]
It follows:
\[
\tr{\bar U \cdot \bar D \cdot \Pi} > \max_{D \in \mathit{\Phi}_{cS}(C)}\tr {U\cdot D \cdot  \mathit{ \Pi}}.
\]
Thus, $C$ cannot be more convexified Shannon-useful than $\bar C$.
%, a contradiction to . Hence, such $\bar D$ cannot exist and therefor $\mathit{\Phi}_{cS}\left(\bar C\right)\subseteq \mathit{\Phi}_{cS}\left(C\right)$ holds.

\underline{\textsc{step 2:}}\; $C\succcurlyeq_{cS} \bar C \Leftarrow \mathit{\Phi}_{cS}(C) \supseteq \mathit{\Phi}_{cS}(\bar C).$\newline
This is clear, since the convexified Shannon-usefulness is defined via maximization over the convexified Shannon policy space.

\underline{\textsc{step 3:}} \;$\mathit{\Phi}_{cS}(C) \supseteq \mathit{\Phi}_{cS}(\bar C) \Rightarrow C \trianglerighteq_{cS} \bar C.$\newline
%Since the convexified Shannon policy space of $\bar C$ is a subset of that of $C$ and
From $ \bar C \in \mathit{\Phi}_{cS}(\bar C)$ it follows by the assumption that $ \bar C \in \mathit{\Phi}_{cS}( C)$. Due to the definition of $\mathit{\Phi}_{cS}( C)$, the channel $\bar C$ is a convexified Shannon-garbling of $C$.

\underline{\textsc{step 4:}} \;$\mathit{\Phi}_{cS}(C) \supseteq \mathit{\Phi}_{cS}(\bar C) \Leftarrow C \trianglerighteq_{cS} \bar C.$\newline
According to the definition of the convexified Shannon-order we have
\[
\bar C= \sum\limits_{j=1}^\ell q_j\cdot M_j \cdot C \cdot N_j,
\]
where $q_i$, $i=1,\ldots, \ell$ is a probability distribution
and $M_i \in \R^{m \times m}$, $N_i \in \R^{n \times n}$, $i=1,\ldots,\ell$,
are stochastic matrices.
We then write for $D \in \mathit{\Phi}_{cS}(\bar C)$:
\[
\begin{array}{rcl}
   D   &=   &\displaystyle \sum\limits_{i=1}^k p_i \cdot A_i \cdot \bar C \cdot B_i = \sum\limits_{i=1}^k p_i \cdot A_i \cdot \left(\sum\limits_{j=1}^\ell q_j\cdot M_j \cdot C \cdot N_j\right)\cdot B_i   \\ \\
   &= &\displaystyle \sum\limits_{i=1}^k \sum\limits_{j=1}^\ell p_i\cdot q_j \cdot A_i \cdot M_j \cdot C \cdot N_j\cdot B_i = \sum\limits_{i=1}^k \sum\limits_{j=1}^\ell \bar p_{ij} \cdot \bar A_{ij} \cdot C \cdot \bar B_{ij} 
\end{array}
\]
with the joint probability distribution
\[
   \bar p_{ij} = p_i\cdot q_j, \quad i=1, \ldots, k, \quad j =1,\ldots, \ell,
\]
and stochastic matrices
\[
   \bar A_{ij} = A_i \cdot M_j, \quad \bar B_{ij} = N_j\cdot B_i,
   \quad i=1, \ldots, k, \quad j =1,\ldots, \ell.
\]
Hence, $D \in \mathit{\Phi}_{cS}(C)$.
\qed

\begin{remark}[Shannon-usefulness]
Without convexification the Shannon policy space for a channel $C$ can be defined as follows:
\[
\mathit{\Phi}_S(C) = \left\{\left.D \in \R^{m \times n}\right\vert \mbox{There exist } A, B \in {\mathcal A} \mbox{ with } D=A \cdot C \cdot B \right\}.
\]
Let $C, \bar C \in \R^{m\times n}$ be two channels with the same input-alphabet.
We say that $C$ is more Shannon-useful than $\bar C$ (denoted by $C \succcurlyeq_S \bar C$) if for all utility matrices $U \in \mathcal U$ it holds:
\[
\max_{D \in \mathit{\Phi}_S(C)}\tr {U\cdot D \cdot  \mathit{ \Pi}}\ge \max_{D \in \mathit{\Phi}_S\left(\bar C\right)}\tr {U\cdot D \cdot  \mathit{ \Pi}}.
\]
Whether Shannon-order can be characterized by this notion of Shannon-usefulness, is not clear. We postpone this question to future research. The main difficulty here is that the Shannon policy space $\mathit{\Phi}_S(C)$ is not convex. Thus, the application of the separation theorem is not possible. \qed
\end{remark}

\bibliographystyle{apalike}
\bibliography{lit.bib}

\end{document}